\documentclass[12pt,leqno]{amsart}

\usepackage{amsmath,amsthm,amsfonts,amscd,amssymb,euscript}
\numberwithin{equation}{section} 
\newcommand{\dbar}{\ensuremath{\bar \partial}}
\newcommand{\dbarb}{\ensuremath{\bar \partial _b }}

\newcommand{\ad}{\ensuremath{\bar \partial _b ^{*}  }}

\newcommand{\adt}{\ensuremath{\bar \partial _{b,t} ^{*}  }}

\newcommand{\adp}{\ensuremath{\bar \partial _b ^{*,\, +}  }}
\newcommand{\adm}{\ensuremath{\bar \partial _b ^{*, \, -}  }}
\newcommand{\ado}{\ensuremath{\bar \partial _b ^{*, \, 0}  }}

\newcommand{\boxbt}{\ensuremath{\square_{b,t}}}
\newcommand{\half}{\ensuremath{\frac {1}{2} \,}}
\newcommand{\C}{\ensuremath{{\mathbb C}}}

\newcommand{\N}{\ensuremath{{\mathbb N}}}

\newcommand{\eps}{\epsilon}

\newcommand{\smooth}{\ensuremath{C^{\infty}}}

\newcommand{\dzi}[1]{\ensuremath{\frac {\partial} {\partial z_{#1}}}}

\newcommand{\dzib}[1]{\ensuremath{\frac {\partial} {\partial \bar z_{#1}}}}

\newcommand{\nul}{\ensuremath{\mathcal{N}}}
\newcommand{\nulperp}{\ensuremath{{\,}^\perp \hspace{-0.02em}\mathcal{N}}}
\newcommand{\ran}{\ensuremath{\mathcal{R}}}

\newcommand{\hol}{\ensuremath{{\mathfrak H}}}
\newcommand{\harmon}{\ensuremath{{\mathcal{H}}_t}}

\newcommand{\harperp}{\ensuremath{{\,}^\perp \hspace{-0.02em}{\mathcal{H}}_t}}
\newcommand{\ze}{\ensuremath{\zeta}}
\newcommand{\dze}{\ensuremath{\tilde{\zeta}}}

\newcommand{\ppt}{\ensuremath{{\Psi}_t^+}}
\newcommand{\pot}{\ensuremath{{\Psi}_t^0}}
\newcommand{\pmt}{\ensuremath{{\Psi}_t^-}}
\newcommand{\pptd}{\ensuremath{{\widetilde{\Psi}}_t^+}}
\newcommand{\potd}{\ensuremath{{\widetilde{\Psi}}_t^0}}

\newcommand{\pmtd}{\ensuremath{{\widetilde{\Psi}}_t^-}}
\newcommand{\pptg}[1]{\ensuremath{{\Psi}_{t, \, #1}^+}}
\newcommand{\potg}[1]{\ensuremath{{\Psi}_{t, \, #1}^0}}
\newcommand{\pmtg}[1]{\ensuremath{{\Psi}_{t, \, #1}^-}}
\newcommand{\pptdg}[1]{\ensuremath{{\widetilde{\Psi}}_{t, \, #1}^+}}
\newcommand{\potdg}[1]{\ensuremath{{\widetilde{\Psi}}_{t, \, #1}^0}}

\newcommand{\pmtdg}[1]{\ensuremath{{\widetilde{\Psi}}_{t, \, #1}^-}}

\newcommand{\sqtnorm}[1]{\ensuremath{{\langle| {#1}|\rangle}^2_t}}
\newcommand{\unsqtnorm}[1]{\ensuremath{{\langle| {#1}|\rangle}_t}}

\newcommand{\tnorm}[2]{\ensuremath{{\langle| {#1},{#2}|\rangle}_t}}

\newcommand{\cp}{\ensuremath{{\mathcal{C}}^+}}
\newcommand{\co}{\ensuremath{{\mathcal{C}}^0}}
\newcommand{\cm}{\ensuremath{{\mathcal{C}}^-}}
\newcommand{\cpd}{\ensuremath{{\widetilde{\mathcal{C}}}^+}}
\newcommand{\cod}{\ensuremath{{\widetilde{\mathcal{C}}}^0}}
\newcommand{\cmd}{\ensuremath{{\widetilde{\mathcal{C}}}^-}}
\newcommand{\qbtg}[1]{\ensuremath{Q_{b,t} ({#1},{#1})}}
\newcommand{\qbtgs}[2]{\ensuremath{Q_{b,t} ({#1},{#2})}}

\newcommand{\qbtp}[1]{\ensuremath{Q^{\, l}_{b,\, +,t} ({#1},{#1})}}
\newcommand{\qbtm}[1]{\ensuremath{Q^{\, l}_{b,\, -,t} ({#1},{#1})}}
\newcommand{\qbto}[1]{\ensuremath{Q^{\, l}_{b,\, 0,t} ({#1},{#1})}}

\begin{document}
\title[The $\dbarb$ equation]{The $\dbarb$ equation on weakly pseudoconvex \\ CR
manifolds of dimension 3}
\author{Joseph J. Kohn and Andreea C. Nicoara}
%\address{Department of Mathematics, Princeton University, Princeton,
%NJ 08544 \\ Department of Mathematics, Harvard University, Cambridge,
%MA 02138}

%\email{kohn@math.princeton.edu; anicoara@math.harvard.edu}

\maketitle
%\begin{abstract}

%Abstract text

%\end{abstract}

\tableofcontents

\section{Introduction} \label{intro}
\medskip
\medskip

The starting point of this work is the second author's thesis
\cite{Andreea} accomplished under the direction of the first
author. In \cite{Andreea} the main result is that on a smooth,
compact, orientable, and weakly pseudoconvex CR manifold of real
dimension $5$ or higher, embedded in $\C^N$ and of codimension one
or higher, the tangential CR operator $\dbarb$ has closed range in
$L^2$ and every Sobolev space $H^s$ with $s>0.$ This is proven
using microlocalization with a specially constructed variant of a
strongly plurisubharmonic function as weight, this function being
the way in which the embeddeding of the manifold is exploited.
Such a microlocalization, however, does not hold for a CR manifold
of real dimension three, whose Levi form is merely a function,
because some positive semi-definite quadratic forms involving the
Levi form that appear in the proof are identically zero in this
case. The aim of this paper is to extend the results of
\cite{Andreea} to the case of a CR manifold of dimension three.

We will employ the same setup as in \cite{Andreea}, which is
described in detail in Section~\ref{defnot}. We will use the same
type of microlocalization, namely the division of the cotangent
space into two truncated cones, a positive and a negative one,
plus a region called the zero portion that overlaps somewhat these
two and on which one obtains the best possible estimates, elliptic
ones. On the positive cone, the estimates for $(0,1)$ forms are
the same as the ones in \cite{Andreea}, but as explained above,
one cannot expect any such estimates to hold on the negative cone.
There, one can only prove estimates for functions. On the zero
region, good estimates hold for both functions and $(0,1)$ forms.
Putting all these estimates together, we are able to prove the
following theorem:

\newtheorem{mainthm}{Theorem}[section]
\begin{mainthm}
Let $M$ be a smooth, compact, orientable, and weakly pseudoconvex
manifold of real dimension $3$ embedded in a complex space $\C^N$
and endowed with the induced CR structure.\label{maintheorem} If
the range of $\dbarb$ is closed in $L^2(M),$ then the range of
$\dbarb$ is also closed in $H^s (M)$ for each $s>0.$ Moreover, for
every $\dbarb$-closed $(0,1)$ form $\alpha \in \smooth$ such that
$\alpha \perp \nul(\ad),$ there exists some function $u \in
\smooth$ such that $\dbarb \, u \, = \, \alpha.$
\end{mainthm}

\noindent Compared to the main result in \cite{Andreea}, this
theorem contains the extra hypothesis that the range of the
tangential CR operator $\dbarb$ is closed in $L^2(M).$ We do not
believe this hypothesis to be necessary; on the contrary, we
conjecture that the embeddability assumption combined with weak
pseudoconvexity should suffice for the range for $\dbarb$ to be
closed in $L^2(M).$ None of the current methods, however, is
strong enough to allow us to prove this conjecture. Until such a
method is developed, we shall justify somewhat this additional
hypothesis by pointing out that Dan Burns proved in
\cite{danburns} that the range of $\dbarb$ fails to be closed in
$L^2$ on the well-known Rossi example of a non-embeddable CR
structure on the three dimensional sphere discussed in
\cite{rossi}.

\bigskip

\bigskip
\section{Definitions and Notation}
\label{defnot}

\medskip
\medskip
\newtheorem{def CR}{Definition}[section]
\begin{def CR}
Let $M$ be a smooth manifold of real dimension $3$ embedded in a
complex space $\C^N.$ An induced CR structure on $M$ is a complex
line bundle of the complexified tangent bundle $\C T(M)$ denoted
by $T^{1,0}(M)$ satisfying $T^{1,0}(M) \, = \, \C T(M) \cap
T^{1,0} (\C^N).$ $M$ endowed with such a CR structure is called a
CR manifold. \label{defCR}
\end{def CR}

\noindent Let $B^q(M)$ be the bundle of $(0,q)$ forms which
consists of skew-symmetric multi-linear maps of $(T^{0,1}(M))^q$
into $\C$. $B^0 (M)$ is then the set of functions on $M$. Since
$M$ is endowed with the induced CR structure coming from $\C^N,$
there exists a natural restriction of the de Rham exterior
derivative to $B^q (M)$, which we will denote by $\dbarb$. Thus,
$\dbarb: \, B^q(M) \rightarrow B^{q+1}(M)$. Clearly, the target
space contains non-zero elements only for $q \, = \, 0.$

\medskip
\noindent Next, it is only natural to choose as a metric the
restriction on $\C T(M)$ of the usual Hermitian inner product on
$\C^N,$ since this Riemannian metric is compatible with the CR
structure on $M,$ namely the spaces $T^{1,0}_P (M)$ and $T^{0,1}_P
(M)$ are orthogonal under it because $\langle \dzi{i}, \dzib{j}
\rangle \, = \, 0$ for all $1 \leq i,j \leq N$. We then define a
Hermitian inner product on $B^q (M)$ by $$(\phi, \psi) = \int
{\langle \phi, \psi \rangle}_z \, dV,$$ where $dV$ denotes the
volume element and $ {\langle \phi, \psi \rangle}_z$ the inner
product induced on $B^q (M)$ by the metric on $\C T(M)$ at each $z
\in M$. Let $||\, \cdot \,||$ be the corresponding norm and $L_2^q
(M)$ the Hilbert space obtained by completing $B^q (M)$ under this
norm. We will only consider the $L^2$ closure of $\dbarb,$ which
we will once again denote it by $\dbarb.$ Its domain is defined as
follows:

\smallskip
\newtheorem{def domdbarb}[def CR]{Definition}
\begin{def domdbarb}
$Dom(\dbarb)$ is the subset of $L_2^q (M)$ composed of all forms
$\phi$ for which there exists a sequence of $\{ \phi_\nu \}_\nu$
in $B^q (M)$ satisfying:
\begin{enumerate}
\item[(i)] $\phi \, = \, \lim_{\nu \rightarrow \infty} \phi_\nu$
in $L^2,$ where $\phi_\nu$ is smooth and
\item[(ii)] $\{ \dbarb \phi_\nu \}_\nu$ is a Cauchy sequence in
$L_2^{q+1} (M)$.
\end{enumerate}
\end{def domdbarb}

\noindent For all $\phi \in Dom(\dbarb)$, let $\lim_{\nu
\rightarrow \infty} \dbarb \phi_\nu \, = \, \dbarb \phi$ which is
thus well-defined. We need to define next $\ad$, the $L^2$ adjoint
of $\dbarb$. Again, we first define its domain:

\smallskip
\newtheorem{def domad}[def CR]{Definition}
\begin{def domad}
$Dom(\ad)$ is the subset of $L_2^q (M)$ composed of all forms
\label{domadjoint} $\phi$ for which there exists a constant $C>0$
such that
$$|(\phi, \dbarb \psi)| \leq C ||\psi||$$ for all $\psi \in
Dom(\dbarb).$
\end{def domad}

\noindent For all $\phi \in Dom(\ad)$, we let $\ad \phi$ be the
unique form in $L_2^q (M)$ satisfying $$(\ad \phi, \psi) = (\phi,
\dbarb \psi),$$ for all $\psi \in Dom(\dbarb)$.

\smallskip
The tangent space to $M$ in the neighborhood $U$, $T(U)$ is
spanned by some $(1,0)$ vector $L$, its conjugate $\overline{L}$,
and one more vector $T$ taken to be purely imaginary, i.e.
$\overline T = -T$. Since we only consider orientable CR
manifolds, the Levi form can be defined globally. We let $\gamma$
be a purely imaginary global $1$-form on $M$ which annihilates
$T^{1,0}(M)\oplus T^{0,1} (M).$  $\gamma$ is not unique, so we
normalize by choosing it in such a way that $\langle \gamma, T
\rangle \, = \, -1.$ Note this implies $\gamma$ is nowhere
vanishing.

\smallskip
\newtheorem{def Levi}[def CR]{Definition}
\begin{def Levi}
The Levi form at a point $z \in M$ is the Hermitian form given by
$\langle d\gamma_z, L \wedge {\overline{L}} \rangle,$ where $L$ is
a vector field in $T_z^{1,0} (U),$ $U$ a neighborhood of $z$ in
$M.$ In dimension 3, the Levi form has only one coefficient, so we
call $M$ weakly pseudo-convex if there exists a form $\gamma$ such
that this coefficient is non-negative at all $z \in M$ and
strongly pseudo-convex if there exists a form $\gamma$ such that
it is positive at all $z \in M.$
\end{def Levi}

\medskip
\noindent If $L$ is the vector field that spans $T^{1,0} (U)$ for
some neighborhoood $U$ of $M,$ let $\omega$ be the $(1,0)$ form
dual to it. In this setup, $\dbarb$ is given by the following: If
$u$ is a smooth function on $U$, then
$$\dbarb (u)= \overline{L}(u) \, \overline{\omega}.$$ If $\varphi \, = \, v \, \overline{\omega}$ is a
$(0,1)$ form on $U$ automatically $\dbarb \, \varphi \, = \, 0.$
As for the $L^2$ adjoint, $\ad$, $$\ad \varphi = \overline{L}^{\,
*} (v) = - L (v) ,$$ where $\overline{L}^{\, *}$ is the $L^2$
adjoint of $\overline{L}.$

We will work with the same inner product on $L_2^q (M)$ as in
\cite{Andreea}, which will be the sum of three inner products: the
inner product without any weight defined above, which we will
denote from now on by ${( \ \cdot \ , \ \cdot \ )}_0$; the inner
product with weight $e^{-t\lambda}$ , ${( \ \cdot \ , \ \cdot \
)}_t \, = \, {(e^{-t\lambda} \ \cdot \ , \ \cdot \ )}_0$; and the
inner product with weight $e^{t\lambda}$, ${( \ \cdot \ , \ \cdot
\ )}_{-t} \, = \, {(e^{t\lambda} \ \cdot \ , \ \cdot \ )}_0$.
Notice that each of these three inner products determines an $L^2$
adjoint, so we denote by $\adp$ the $L^2$ adjoint on ${( \ \cdot \
, \ \cdot \ )}_t$, by $\adm$ the $L^2$ adjoint on ${( \ \cdot \ ,
\ \cdot \ )}_{-t}$, and by $\ado$ the $L^2$ adjoint on ${( \ \cdot
\ , \ \cdot \ )}_0$, although since there is no weight function in
this case, $\ado$ equals precisely $\ad$. Now let us compute
$\adp$ and $\adm$. Let $u$ be a smooth function and $\varphi \, =
\,v \, \overline{\omega}$ a smooth $(0,1)$ form, then
$$(\dbarb \, u, \varphi)_t = (\dbarb \, u, e^{-t \lambda} \varphi)_0 =
(u, \ad (e^{-t \lambda} \varphi))_0 = (u,e^{-t \lambda}(\ad -t
[\ad, \lambda])\varphi)_0 = (u,( \ad -t [\ad,\lambda])\varphi)_t
\, . $$ Therefore, $\adp \varphi =(\ad -t [\ad,\lambda]) \varphi =
\overline{L}^{\, *, \, t} (\varphi_i) =  -  L (v) + t L
(\lambda)\, v ,$ where $\overline{L}^{\, *, \, t}$ is the $L^2$
adjoint of $\overline{L}$ with respect to ${( \ \cdot \ , \ \cdot
\ )}_t$. Similarly,$$(\dbarb \, u, \varphi)_{-t} = (\dbarb \, u,
e^{t \lambda} \varphi)_0 = (u, \ad (e^{t \lambda} \varphi))_0 =
(u,e^{t \lambda}( \ad + t[\ad, \lambda])\varphi)_0 = (u,( \ad +
t[\ad,\lambda])\varphi)_{-t} \, .
$$ So then $\adm \varphi = (\ad + t[\ad,\lambda]) \varphi =
\overline{L}^{\, *, \, -t} (v) =  -  L (v) - t L (\lambda)\, v, $
where $\overline{L}^{\, *, \, -t}$ is the $L^2$ adjoint of
$\overline{L}$ with respect to ${( \ \cdot \ , \ \cdot \ )}_{-t}$.
Just as in \cite{Andreea}, $\lambda$ is chosen to be CR
plurisubharmonic, which is defined as follows:

\smallskip
\newtheorem{def CRpshglobal}[def CR]{Definition}
\begin{def CRpshglobal}
Let $M$ be a CR manifold. A $\smooth$ real-valued function
$\lambda$ defined in the \label{CRpshdefglobal} neighborhood of
$M$ is called strongly CR plurisubharmonic if $\, \exists$ $A_0>0$
such that $\langle  \frac{1}{2} \, (\partial_b \dbarb \lambda -
\dbarb \partial_b \lambda)+ A_0 \,d\gamma, L \wedge {\overline{L}}
\rangle$ is strictly positive $\, \forall$ $L \in T^{1,0} (M),$
where $\langle d\gamma, L \wedge {\overline{L}} \rangle$ is the
invariant expression of the Levi form. $\lambda$ is called weakly
CR plurisubharmonic if $\, \langle \frac{1}{2} \, (\partial_b
\dbarb \lambda - \dbarb \partial_b \lambda)+ A_0 \, d\gamma, L
\wedge {\overline{L}} \rangle$ is just non-negative.
\end{def CRpshglobal}

\smallskip\noindent {\bf Remark.} Note that the previous definition is
trivially satisfied for any $\lambda$ if $M$ is strongly
pseudoconvex. This reflects the fact that in such a case the
microlocal argument can be carried out in absence of a weight
function.

\smallskip
\noindent The following properties of CR plurisubharmonic
functions stated here for 3-dimensional CR manifolds were proven
for $(2n-1)$-dimensional CR manifolds in \cite{Andreea} and are
relevant for the upcoming argument:

\smallskip
\newtheorem{CRpshresultsummary}[def CR]{Proposition}
\begin{CRpshresultsummary}
Let $M$ be a compact, smooth, orientable, \label{CRpshprops}
weakly pseudoconvex CR manifold of real dimension 3 embedded in a
complex space $\C^N$ and endowed with an induced CR structure.
\begin{enumerate}
\item[(i)] Around each point $P \in M$, there exists a
small enough neighborhood $U$ and a local orthonormal basis $L,
L',\overline{L},\overline{L}'$ of the $2$ dimensional complex
bundle containing $T M$ when restricted to $U$, satisfying $[L,
\overline{L}] \big|_P \, = \, c \,T$, where $T \, = \, L' -
\overline{L}'$ and $c$ is the coefficient of the Levi form in the
local basis $L, \overline{L}, T$ of $T M$.
\item[(ii)] If $\lambda$ defined on $M$ is strongly CR plurisubharmonic, then $\, \exists$
$A_0>0$ such that $ \frac {1}{2} \, ( \overline{L} L (\lambda)+ L
\overline{L} (\lambda)) +A_0 \, c$ is positive in a neighborhood
$U'$ around $P$ which is smaller than $U$. $A_0$ is of course
independent of $P$ or $U$, and the size of $U'$ depends on it. If
$\lambda$ is weakly CR plurisubharmonic, then $ \frac {1}{2} \, (
\overline{L} L (\lambda)+ L \overline{L} (\lambda)) +A_0 \, c$ is
just non-negative.
\item[(iii)] If $\lambda$ is strongly
 plurisubharmonic on $\C^N$, then $\lambda$ is also strongly CR plurisubharmonic on $M$.
\item[(iv)]  Let $$\overline{L}^{\,*, \, \pm t} = - L  \pm t L
(\lambda).$$ For each small positive number $\eps_G,$ there exists
a covering $\{ V_\mu \}_\mu,$ a local basis $L, \overline{L}, T$
of $T V_\mu$ for each $\mu,$ and $\smooth$ functions $a,$ $b,$
$g,$ and $e$ so that the bracket $[L,\overline{L}]$ has the
following form:
\begin{equation}
\begin{split}
[L,\overline{L}] &= c T + a \overline{L} + b \overline{L}^{\,*, \,
\pm t}\pm t \, e \label{bracket},
\end{split}
\end{equation}
where $|e|$ is bounded independently of $t,$ namely $$|e| \leq
\eps_G.$$
\end{enumerate}
\end{CRpshresultsummary}

\noindent Just as in \cite{Andreea}, we will use microlocalization
to prove the main estimates, so the following definitions are
necessary:

\smallskip
\newtheorem{def dom}[def CR]{Definition}
\begin{def dom}
Let $P$ be a pseudodifferential operator of order zero, then
another pseudodifferential operator of order zero, $\tilde{P}$ is
said to dominate $P$, if the symbol of $\tilde{P}$ is identically
equal to $1$ on a neighborhood of the support of the symbol of $P$
and the support of the symbol of $\tilde{P}$ is slightly larger
than the support of the symbol of $P$.
\end{def dom}

\noindent In particular, the previous definition also applies to
cutoff functions, all of which are pseudodifferential operators of
order zero.

\smallskip
\newtheorem{def tdep}[def CR]{Definition}
\begin{def tdep}
Let $P(t)$ be a family of pseudodifferential operators depending upon
a parameter $t$. Such a family is called zero order $t$ dependent if
it can be written as $P_1+t \, P_2$, where $P_1$ and $P_2$ are
pseudodifferential operators independent of $t$ and $P_2$ has order zero.
\end{def tdep}

\smallskip
\newtheorem{def invtdep}[def CR]{Definition}
\begin{def invtdep}
Let $P(t)$ be a family of pseudodifferential operators depending upon
a parameter $t$ such that $P(t)$ has order zero. Such a family $P(t)$ is called inverse
zero order $t$ dependent if its symbol $\sigma(P)$ satisfies
$D^{\alpha}_{\xi} \sigma (P) \, = \, D^{\alpha}_{\xi} p(x, \xi) \,
= \, \frac {1}{t^{|\alpha|}} \, q(x,\xi)$ for $|\alpha| \geq 0$,
where $q(x,\xi)$ is bounded independently of $t$.
\end{def invtdep}

\noindent Notice that the two adjoint operators $\adp$ and $\adm$
are zero order $t$ dependent, according to the definition given
above.

\bigskip
\noindent The microlocalization consists in dividing the Fourier
transform space into three conveniently chosen regions, two
truncated cones $\cp$ and $\cm$ and another region $\co$, with
some overlap. Let the coordinates on the Fourier transform space
be $\xi \, = \, (\xi_1, \xi_2, \xi_3)$. Write $\xi' \, = \,(\xi_1,
\xi_2)$, so then $\xi \, = \, (\xi', \xi_3)$. The work is done in
coordinate patches on $M,$ each of which has defined on it local
coordinates such that $\xi'$ is dual to the holomorphic part of
the tangent bundle $T^{1,0}(M)\oplus T^{0,1} (M)$ and $\xi_3$ is
dual to the totally real part of the tangent bundle of $M$ spanned
by the "bad direction," $T.$ Duality in this context merely means
that the pairing between the specified cotangent and tangent
sub-bundles is nondegenerate at each point of the coordinate
patch. Define $\cp \, = \, \{ \xi \, | \, \xi_3 \geq \frac {1}{2}
|\xi'| \: and \: |\xi| \geq \ 1 \}$. Then $\cm \, = \, \{ \xi \, |
\, -\xi \in \cp \}$, and finally $\co \, = \, \{ \xi \, | \,
-\frac {3}{4} |\xi'| \leq \xi_3 \leq \frac {3}{4} |\xi'| \} \cup
\{ \xi \, | \, |\xi| \leq 1 \}$. Notice that by definition, $\cp$
and $\co$ overlap on two smaller cones and part of the sphere of
radius $1$ and similarly $\cm$ and $\co$, whereas $\cp$ and $\cm$
do not intersect.

Let us now define three functions on $\{ |\xi'|^2 + |\xi_3|^2 \, =
\, 1 \}$, which is the unit sphere in $\xi$ space. $\psi^+$,
$\psi^-$, and $\psi^0$ are smooth, take values in $[0,1]$, and
satisfy the condition of symbols of pseudodifferential operators
of order zero. Moreover, $\psi^+$ is supported in $\{ \xi \,  | \,
\xi_3 \geq \frac {1}{2} |\xi'| \}$ and $\psi^+ \, \equiv \, 1$ on
the subset $\{ \xi \, | \, \xi _3 \geq \frac {3}{4} |\xi'| \}$.
Then let $\psi^- (x, \xi)\, = \, \psi^+ (x, -\xi)$ which means
that $\psi^-$ is supported in $\{ \xi \,  | \, \xi_3 \leq - \frac
{1}{2} |\xi'| \}$ and $\psi^- \, \equiv \, 1$ on the subset $\{
\xi \, | \, \xi _3 \leq -\frac {3}{4} |\xi'| \}$. Finally, let
$\psi^0 (\xi)$ satisfy $(\psi^0 (\xi))^2 \, = \, 1 - (\psi^+
(\xi))^2 - (\psi^- (\xi))^2$ which means that $\psi^0$ is
supported in $\co$ and $\psi^0 \, \equiv \, 1$ on the subset $\{
\xi \, | \, -\frac {1}{2} |\xi'| \leq \xi_3 \leq \frac {1}{2}
|\xi'| \}$. Next extend $\psi^+$, $\psi^-$, and $\psi^0$
homogeneously by setting $\psi^+ (\xi) \, = \, \psi^+ (\frac
{\xi}{|\xi|})$, for $\xi$ in $\cp$ outside of the unit sphere.
Similarly, let $\psi^- (\xi) \, = \, \psi^- (\frac {\xi}{|\xi|})$,
for $\xi$ in $\cm$ outside of the unit sphere and $\psi^0 (\xi) \,
= \, \psi^0 (\frac {\xi}{|\xi|})$, for $\xi$ in $\co$ outside of
the unit sphere. Extend $\psi^+$, $\psi^-$, and $\psi^0$ inside
the unit sphere in some smooth way so that $(\psi^+)^2 + (\psi^-
)^2+ (\psi^0)^2 \, = \, 1$ still holds. Now we have the functions
$\psi^+$, $\psi^-$, and $\psi^0$ defined everywhere on $\xi$
space. We then define $\psi^+_t (\xi) \, = \, \psi^+ (\frac
{\xi}{tA})$, $\psi^0_t (\xi) \, = \, \psi^0 (\frac {\xi}{tA})$,
and $\psi^-_t (\xi) \, = \, \psi^- (\frac {\xi}{tA})$ for some
positive constant $A$ to be chosen later. Let $\ppt,$ $\pot,$ and
$\pmt$ be pseudodifferential operators of order zero with symbols
$\psi^+_t$, $\psi^-_t$, and $\psi^0_t$ respectively. By
construction, the following is true: $$(\ppt)^* \ppt + (\pot)^*
\pot  +(\pmt)^* \pmt = Id,$$ modulo a smoothing operator. The
global norm is defined using a very special covering of $M,$ $\{
U_\nu \}_\nu,$ satisfying a number of technical conditions proven
in \cite{Andreea} and summarized in the next lemma:

\smallskip
\newtheorem{coveringlemma}[def CR]{Lemma}
\begin{coveringlemma}
Let $M$ be a compact, smooth, orientable, weakly pseudoconvex CR
manifold of real dimension 3 embedded in a complex space $\C^N$
and endowed with an induced CR structure. \label{coveringprops}
Let $\pptg{\nu}$, $\potg{\nu}$, and $\pmtg{\nu}$ be the
pseudodifferential operators of order zero defined on $U_\nu$ for
each $\nu$ and $\cp_\nu$, $\co_\nu$, and $\cp_\nu$ be the three
regions of the $\xi$ space dual to $U_\nu$ on which the symbol of
each of those three pseudodifferential operators is supported.
Moreover, define $\pptdg{\nu}$ and $\pmtdg{\nu}$ so that they
dominate $\pptg{\nu}$ and $\pmtg{\nu}$ respectively and are also
inverse zero order $t$ dependent. We denote by $\cpd_\nu$ and
$\cmd_\nu$ the supports of the symbols of $\pptdg{\nu}$ and
$\pmtdg{\nu}$, respectively.
\begin{enumerate}
\item[(i)] If $\{ V_\mu \}_\mu$ is the covering of $M$ given by
part (iv) of Proposition~\ref{CRpshprops}, then $\{ U_\nu \}_\nu$
is such that for each $\nu,$ $\exists \: \mu(\nu)$ with the
property that $U_\nu$ and all $U_\eta$ satisfying $U_\nu \cap
U_\eta \, \neq \, \emptyset$ are contained in the neighborhood
$V_{\mu(\nu)}$ of $\{ V_\mu \}_\mu;$
\item[(ii)] Let $U_\nu$ and $U_\mu$ be two neighborhoods such that
$U_\nu \cap \, U_\mu \, \neq \, \emptyset$. There exists a
diffeomorphism $\vartheta$ between $U_\nu$ and $U_\mu$ with
Jacobian $\mathcal{J}_\vartheta$ satisfying ${\,}^t
\hspace{-0.25em} \mathcal{J}_\vartheta (\cp_\mu) \cap \cm_\nu \, =
\, \emptyset$ and $\cp_\nu \cap \, {\,}^t \hspace{-0.25em}
\mathcal{J}_\vartheta (\cm_\mu) \, = \, \emptyset$, where ${{\,}^t
\hspace{-0.25em} \mathcal{J}_\vartheta}$ is the inverse of the
transpose of the Jacobian of $\vartheta$;
\item[(iii)] Let ${\,}^\vartheta \pptg{\mu}$, ${\,}^\vartheta \pmtg{\mu}$,
and ${\,}^\vartheta \potg{\mu}$ be the transfers of $\pptg{\mu}$,
$\pmtg{\mu}$, and $\potg{\mu}$ respectively via $\vartheta$, then
on $\{ \xi \, | \, \xi_{2n-1} \, \geq \, \frac {4}{5} |\xi'| \:
and \: |\xi| \geq (1 + \eps)  \, tA \}$, the principal symbol of
${\,}^\vartheta \pptg{\mu}$ is identically equal to  $1$, on $\{
\xi \, | \, \xi_{2n-1} \, \leq \, - \frac {4}{5} |\xi'| \: and \:
|\xi| \geq (1+ \eps) \, tA \}$, the principal symbol of
${\,}^\vartheta \pmtg{\mu}$ is identically equal to $1$, and on
$\{ \xi \, | \, - \frac {1}{3} |\xi'| \, \leq \, \xi_{2n-1} \,
\leq \, \frac {1}{3} |\xi'| \: and \: |\xi| \geq (1 + \eps) \, tA
\}$, the principal symbol of ${\,}^\vartheta \potg{\mu}$ is
identically equal to $1$, where $\eps > 0$ and can be very small;
\item[(iv)] ${\,}^t \hspace{-0.25em}
\mathcal{J}_\vartheta (\cpd_\mu) \cap \cm_\nu \, = \, \emptyset$
and $\cp_\nu \cap \, {\,}^t \hspace{-0.25em} \mathcal{J}_\vartheta
(\cmd_\mu) \, = \, \emptyset$;
\item[(v)] Let ${\,}^\vartheta \pptdg{\mu}$, ${\,}^\vartheta
\pmtdg{\mu}$ be the transfers via $\vartheta$ of $\pptdg{\mu}$ and
$\pmtdg{\mu}$, respectively. Then the principal symbol of
${\,}^\vartheta \pptdg{\mu}$ is identically $1$ on $\cp_\nu$ and
the principal symbol of ${\,}^\vartheta \pmtdg{\mu}$ is
identically $1$ on $\cm_\nu$;
\item[(vi)] $\cpd_\mu \cap \cmd_\mu \, = \, \emptyset$.
\end{enumerate}
\end{coveringlemma}

\smallskip
\noindent Henceforth, the left superscript $\vartheta$ indicating
the transfer of a pseudodifferential operator into another local
coordinate system will be suppressed to simplify our notation. Let
now $\{ \ze_\nu \}_\nu$ be a partition of unity subordinate to the
covering $\{ U_\nu \}_\nu$ satisfying $\sum_\nu {\ze_\nu}^2 \, =
\, 1$, and for each $\nu$ let $\dze_\nu$ be a cutoff function that
dominates $\ze_\nu$ such that $supp \, (\dze_\nu) \subset U_\nu$.
Then we define the global norm as follows:
$${\langle| \varphi |\rangle}_t^2 =
\sum_\nu \: (||\dze_\nu \pptg {\nu} \ze_\nu v^\nu||_t^2 +
||\dze_\nu \potg {\nu} \ze_\nu v^\nu ||_0^2 + ||\dze_\nu \pmtg
{\nu} \ze_\nu v^\nu||_{-t}^2)$$ where $\varphi \, = \, v \,
\overline{\omega} $ is a $(0,1)$ form in $L_2^1 (M)$ and $v^\nu$
is the coefficient of the form expressed in the local coordinates
on $U_\nu$. We define the norm exactly the same way for functions.
The following two facts are proven in \cite{Andreea}:

\smallskip
\newtheorem{normpropslemma}[def CR]{Lemma}
\begin{normpropslemma}
Let $M$ be a compact, smooth, orientable, \label{normprops} weakly
pseudoconvex CR manifold of real dimension 3 embedded in a complex
space $\C^N$ and endowed with an induced CR structure.
\begin{enumerate}
\item[(i)] For any $t$, there exist two positive
constants depending on $t$, $C_t$ and $C'_t$ such that
$$ C_t ||\varphi||_0^2 \leq {\langle| \varphi |\rangle}_t^2 \leq
C'_t  ||\varphi||_0^2,$$ where $\varphi$ is a form in $L_2^1 (M)$;
\item[(ii)] There exists a self-adjoint operator $G_t$ such
that $$ (\varphi, \phi)_0 = {\langle| \varphi, G_t \, \phi
|\rangle}_t,$$ for any two forms $\varphi$ and $\phi$ in $L_2^1
(M)$. $G_t$ is the inverse of the operator $F_t$ given by
$$\sum_\nu \: \left(\ze_\nu (\pptg {\nu})^* \dze_\nu e^{-t
\lambda} \dze_\nu \pptg {\nu} \ze_\nu + \ze_\nu (\potg {\nu})^*
\dze_\nu^2 \potg {\nu} \ze_\nu + \ze_\nu (\pmtg {\nu})^* \dze_\nu
e^{t \lambda} \dze_\nu \pmtg {\nu} \ze_\nu\right).$$
\end{enumerate}
\end{normpropslemma}

\smallskip
\noindent Finally, here is the definition of the Sobolev norm for
forms in this context. Clearly, we make the same definition for
functions.

\smallskip
\newtheorem{Sobolevdef}[def CR]{Definition}
\begin{Sobolevdef} Let the Sobolev norm of order $s$ for a form $\varphi \, = \, v \, \overline{\omega}$
supported on $M$ be given by:
$$||\varphi||^2_s = \sum_\eta  ||
\dze_\eta \, \Lambda^s \, \ze_\eta \, v^\eta||^2_0,$$ where as
usual $\Lambda$ is defined to be the pseudodifferential operator
with symbol $(1+|\xi|^2)^\half.$ Then $$H^s = \{ \varphi \in
B^1(M) \: \big| \: ||\varphi||_s < + \infty \}.$$
\end{Sobolevdef}

\smallskip \noindent A few more definitions are needed to lay out
the terminology used in the statement of
Theorem~\ref{maintheorem}.

\bigskip
\newtheorem{closerangedef}[def CR]{Definition}
\begin{closerangedef}
An operator $P$ has closed range if $\forall \: \alpha \in
\overline{\ran(P)},$ where $\overline{\ran(P)}$ is the closure of
the range of $P,$ $\alpha \in \ran(P).$
\end{closerangedef}

\smallskip \noindent In particular, for $\dbarb$ to have closed
range in $L^2(M)$ means that if we denote by $L^1_2(M)$ the set of
$(0,1)$ forms in $L^2$ of $M,$ the set of closed, $(0,1)$ forms in
$L^2$ decomposes as follows:
$$L^1_2(M) \cap \nul (\dbarb) = \ran(\dbarb) \oplus \harmon,$$
where for each $t$ the corresponding harmonic space is
$$\harmon = \{ \varphi \in Dom(\dbarb) \cap Dom(\ad) \; \big| \;
\dbarb \, \varphi = 0 \: and \: \adt \, \varphi=0 \}.$$ In other
words, given some $\dbarb$-closed $(0,1)$ form $\alpha$ that is in
$L^2$ of $M$ and such that $\alpha \perp \harmon,$ where
orthogonality is defined with respect to the ${\langle| \, \cdot
\, |\rangle}_t$ norm, there exists some function $u_t$ in $L^2$ of
$M$ such that $$\dbarb u_t = \alpha.$$ Note that we index the
function by $t$ since the norm depends on $t,$ hence such a
function exists for each $t,$ and functions corresponding to
different values of $t$ are not necessarily equal.

The closed range property is equivalent to three other properties
because $\dbarb$ is a linear, closed, densely defined operator on
a Hilbert space. See for example \cite{Hormander}:

\smallskip
\newtheorem{closeranequiv}[def CR]{Theorem}
\begin{closeranequiv}
The following four conditions \label{closedrange} are equivalent:
\begin{enumerate}
\item[(i)] $\ran(\dbarb)$ is closed in $H^s$;
\item[(ii)] There exists a $t$-dependent constant $C_t$ such that
\begin{equation}
{\langle| \,\Lambda^s u_t \, |\rangle}_t \leq  C_t \, {\langle| \,\Lambda^s \alpha
\, |\rangle}_t, \label{dbarbrangeineq}
\end{equation}
where $u_t \in Dom(\dbarb)\cap \overline{\ran(\adt)} \subset
Dom(\dbarb) \cap \nulperp(\dbarb)$ and $\dbarb \, u_t \, = \,
\alpha;$
\item[(iii)] $\ran(\adt)$ is closed in $H^s$;
\item[(iv)]  There exists a $t$-dependent constant $C_t$ such that
\begin{equation}
{\langle| \,\Lambda^s \varphi_t \, |\rangle}_t \leq  C_t \, {\langle|
  \, \Lambda^s
u_t \, |\rangle}_t, \label{adtrangeineq}
\end{equation}
where $\varphi_t \in Dom(\adt)\cap \overline{\ran(\dbarb)}\subset
Dom(\adt)\cap\nulperp(\adt) $ and $\adt \, \varphi_t \, = \, u_t.$
\end{enumerate}
The best constants in \ref{dbarbrangeineq} and in
\ref{adtrangeineq} are the same.
\end{closeranequiv}

\smallskip
\noindent In particular, it follows that the $\adt$ problem can be
solved, i.e. for every function $u_t \perp \harmon,$ there exists
a $(0,1)$ form $\varphi_t$ in $L^2,$ $\varphi_t \perp \nul
(\adt),$ such that $$\adt\, \varphi_t = u_t.$$ On $\harperp$ we
define the Kohn Laplacian
$$\boxbt = \dbarb \, \adt + \adt \, \dbarb.$$ $\boxbt$ is
self-adjoint, and on $(0,1)$ forms on a manifold $M$ of dimension
$3,$ $\boxbt \, = \, \dbarb \, \adt.$ Moreover, by the
considerations above, for every $(0,1)$ form $\alpha \perp
\harmon,$ there exists a $(0,1)$ form $\varphi_t \perp \harmon$
such that $$\boxbt \, \varphi_t=\alpha.$$ Let $N_t$ be the
solution operator which takes $\alpha$ to $\varphi_t$. Define
$N_t$ to be identically zero on $\harmon.$ \ref{dbarbrangeineq}
and \ref{adtrangeineq} together imply that $${\langle| \,
\varphi_t \, |\rangle}_t = {\langle| \, N_t \, \alpha \,
|\rangle}_t \leq C'_t \, {\langle| \, \alpha \, |\rangle}_t,$$ for
some $t$-dependent constant $C'_t,$ i.e. the solution operator
$N_t$ is bounded on $L^2.$ Furthermore, notice that $u_t \, = \,
\adt \, N_t \, \alpha$ solves the $\dbarb$ problem $$\dbarb \, u_t
= \alpha.$$

\smallskip\noindent Let us now relate $\nulperp(\adt)$ to
$\nulperp(\ad)$ in order to simplify the hypotheses of our main
result, Theorems~\ref{maintheorem}.

\smallskip
\newtheorem{nulperplemma}[def CR]{Lemma}
\begin{nulperplemma}
Let $\perp$ denote \label{nulperp} perpendicularity with respect
to the unweighted norm $\| \, \cdot \, \|_0,$ and let $\perp_t$
denote perpendicularity with respect to the weighted norm
$\unsqtnorm{ \, \cdot \,}.$ $\alpha \perp \nul(\ad)$ implies
$\alpha \perp_t \nul(\adt)$ for each $t$ and each form or function
$\alpha.$
\end{nulperplemma}

\smallskip
\noindent {\bf Proof:} For any $\psi \in \nul(\adt),$
$$\tnorm{\alpha}{\psi}=\|\alpha, F_t \, \psi \|_0.$$ Showing that $F_t \, \psi \in \nul(\ad)$
would conclude the proof of the lemma. By Equation $4.7$ of
\cite{Andreea}, $\adt$ and $\ad$ are related via the operators
$F_t$ and $G_t$ from Lemma~\ref{normprops} as follows: $$\adt =
\ad + G_t \, [\ad, F_t].$$ This means $$\adt \, \psi = 0 = \ad \,
\psi + G_t \, \ad \, F_t \, \psi  -\ad \, \psi,$$ i.e. $\ad \, F_t
\, \psi \, = \, 0$ a.e. because the operators $F_t$ and $G_t$ give
the correspondence between two equivalent norms, $\| \, \cdot \,
\|_0$ and $\unsqtnorm{ \, \cdot \,},$ so their kernels consist
only of those functions and forms that are zero almost everywhere.
\qed

\smallskip

\smallskip
\newtheorem{Szegoproj}[def CR]{Definition}
\begin{Szegoproj}
Let $\hol(M) = \{f \in Dom(\dbarb) \: \big| \: \dbarb f = 0\}.$
Then the weighted Szeg\"{o} projection $S_{b,t}$ is the operator
that projects $L^2(M)$ to $\hol(M)$ in the ${\langle | \, \cdot \,
|\rangle}_t$ norm.
\end{Szegoproj}

\smallskip\noindent Just as in \cite{Kohnsurvey}, the weighted
Szeg\"{o} projection can be expressed as follows:

\smallskip
\newtheorem{Szegolemma}[def CR]{Lemma}
\begin{Szegolemma}
\label{Szegoexpr} Let $M$ be a compact, orientable, weakly
pseudoconvex CR manifold of dimension $3$ embedded in a complex
space $\C^N$ and endowed with the induced CR structure. If the
range of $\dbarb$ is closed in $L^2(M),$ then
$$S_{b,t} = Id - \adt \, N_t \, \dbarb.,$$ where $Id$ is the identity
operator.
\end{Szegolemma}

\smallskip
\noindent {\bf Proof:} We consider the two complementary cases, $f
\in \hol(M)$ and $f \perp \hol(M),$ and prove this expression for
each. First, if $f \in \hol(M),$ then $(Id - \adt \, N_t \,
\dbarb)\,f \, = \, f$ as expected. If $f \perp \hol(M),$ then
$\adt \, N_t \, \dbarb \, f \, = \, f$ since the ranges of
$\dbarb$ and $\adt$ are closed. Thus, $S_{b,t} \, f \, = \, 0.$
\qed

\medskip\noindent We want to investigate the regularity of the operator
$S_{b,t}$.

\bigskip

\section{Regularity of the Weighted Szeg\"{o} Projection}
\label{szegosection}

\bigskip
We start by recalling two results from \cite{Andreea}, which we
state in the notation appropriate for dimension three. These
results are both local, namely we restrict the discussion to a
small neighborhood $U$ in which all the necessary properties hold,
namely part (iv) of Proposition~\ref{CRpshprops}. Let $\ze$ be a
cutoff function with support in $U$ and $\dze$ a cutoff function
which dominates it, whose support is contained in a slightly
larger open set $U'$.  $U'$ is small enough to allow the existence
of the three pseudodifferential operators of order zero $\ppt$,
$\pot$, and $\pmt$ supported in $\cp$, $\co$, and $\cm$
respectively, as detailed in the previous section. Moreover, let
$\pptd,$ $\potd,$ and $\pmtd$ be pseudodifferential operators of order
zero whose symbols dominate those of  $\ppt$,
$\pot$, and $\pmt$ and are supported in $\cpd$, $\cod$, and $\cmd$
respectively, which are slightly larger than $\cp,$ $\co,$ and $\cm.$ We shall work
with the energy form $$\qbtg{\phi} = \sqtnorm{\dbarb \phi} +
\sqtnorm{\adt \, \phi}$$ and its equivalents for the $t$ norm and
the $-t$ norm,
$$\qbtp {\phi} = {||\dbarb \phi||}_{\, t}^2 + {||\adp \phi||}_{\,
t}^2,$$ and
$$\qbtm {\phi} = {||\dbarb \phi||}_{-t}^2 + {||\adm
\phi||}_{-t}^2.$$

\smallskip
\newtheorem{qbtpcalclemma}{Lemma}[section]
\begin{qbtpcalclemma}
Let $\varphi  \in Dom(\dbarb) \cap Dom(\ad)$ be a $(0,1)$ form
supported in $U',$ \label{qbtp} a neighborhood of a compact,
three-dimensional, weakly-pseudoconvex CR-manifold $M.$ Let there
exist a basis $L, \overline{L}, T$ of $T (U)$ such that $[L,
\overline{L}]= c \, T + a\overline{L} + b \overline{L}^{\, *, \, t} +t
\, e,$ where $a$, $b$, and $e$ are \smooth functions independent
of $t$ and $c$ is the coefficient of the Levi form. Let $M$ be
also endowed with a strongly CR plurisubharmonic function
$\lambda$. Then there exists a constant $C$ independent of $t,$ a $t$-dependent constant $C_t$ and a
positive number $T_0$ such that for any $t \geq T_0$
\begin{equation*}
\begin{split}
\qbtp{\dze \ppt \varphi}+
C_t \, ||\dze \potd \varphi||^2_{\,0}  &\geq C\,t \,||\dze
\ppt\varphi||^2_{\, t}
\end{split}
\end{equation*}
\end{qbtpcalclemma}

\smallskip
\newtheorem{qbtmgardinglemma}[qbtpcalclemma]{Lemma}
\begin{qbtmgardinglemma}
Let $u$ be a function supported in $U'$ such that up to a
smooth term its Fourier transform $\hat u$ is supported in
$\cm,$ and let $g$ be a non-negative function, then the following holds
\label{mgarding}:$$\Re\{(g(-T) u, u )_{-t}\} \geq t \, A \,
(g \, u, u)_{-t}+O(||u||^2_{-t})+O_t(||\dze \potd
u||^2_{\, 0})$$
\end{qbtmgardinglemma}

\smallskip
\noindent The crucial difference between the three-dimensional CR manifold case
and the case of manifolds of dimension five and above is in the handling
of the microlocalization on $\cm,$ so we shall now obtain an
estimate for $\qbtm{\, \cdot\,}$ for functions:

\smallskip
\newtheorem{qbtmcalclemma}[qbtpcalclemma]{Lemma}
\begin{qbtmcalclemma}
Let $u  \in Dom(\dbarb) \cap Dom(\ad)$ be a function
supported in $U',$ \label{qbtm} a neighborhood of a compact,
three-dimensional, weakly-pseudoconvex CR-manifold $M.$ Let there
exist a basis $L, \overline{L}, T$ of $T (U)$ such that $[L,
\overline{L}]= c \, T + a\overline{L} + b \overline{L}^{\, *, -t} -t
\, e,$ where $a$, $b$, and $e$ are \smooth functions independent
of $t$ and $c$ is the coefficient of the Levi form. Let $M$ be
also endowed with a strongly CR plurisubharmonic function
$\lambda$. Then there exists a constant $C'$ independent of $t,$ a $t$-dependent constant $C'_t$ and a
positive number $T'_0$ such that for any $t \geq T'_0$
\begin{equation*}
\begin{split}
\qbtm{\dze \pmt u}+
C'_t \, ||\dze \potd u||^2_{\,0}  &\geq C' \,t \,||\dze
\pmt u||^2_{\,-t}
\end{split}
\end{equation*}
\end{qbtmcalclemma}

\smallskip \noindent {\bf Proof:} Define the pseudodifferential
operators $\ppt,$ $\pot,$ and $\pmt$ using a constant $A \geq A_0,$
where $A_0$ is the CR plurisubharmonicity constant of $\lambda.$
\begin{equation*}
\begin{split}
\qbtm{u} &= ||\dbarb \, u||^2_{-t} = (\overline{L} u, \overline{L}
u)_{-t} = (\overline{L}^{\, *, -t} \overline{L} u, u)_{-t} =
||\overline{L}^{\, *, -t} u||^2_{-t} + ([\overline{L}^{\, *, -t},
  \overline{L}]u, u)_{-t}\\&=||\overline{L}^{\, *, -t} u||^2_{-t} +
([-L-t\, L(\lambda), \overline{L}]u, u)_{-t} \geq (1-\eps) \,
||\overline{L}^{\, *, -t} u||^2_{-t}-\Re\{(c \, T u,
u)_{-t}\}\\&\quad+O(||u||^2_{-t})+t \, \Re\{(\overline{L}L(\lambda) u,
u)_{-t}+t \, \Re\{(e \, u, u)_{-t}\},
\end{split}
\end{equation*}
for some $1 \gg \eps >0.$ Now, replace $u$ by $\dze \pmt u$ in the
previous expression, which has the property that its Fourier transform
is supported in $\cm$ up to a smooth error term. It is now possible to
apply the previous lemma to $-\Re\{(c \, T u, u)_{-t}\}$ to conclude that
\begin{equation*}
\begin{split}
\qbtm{\dze \pmt u} &\geq (1-\eps) \,
||\overline{L}^{\, *, -t}\dze \pmt  u||^2_{-t}+ t \, A_0 (c \, \dze \pmt u,
 \dze \pmt  u)_{-t}+O(||\dze \pmt  u||^2_{-t})\\&\quad+O_t(||\dze \potd  u||^2_{-t}) +t \,
 \Re\{(\overline{L}L(\lambda)\dze \pmt  u,
\dze \pmt  u)_{-t}\}\\&\quad+t \, \Re\{(e \, \dze \pmt u,\dze \pmt u)_{-t}\}.
\end{split}
\end{equation*}
Since the inner product is Hermitian, it is easily seen that
$$ \Re\{(\overline{L}L(\lambda)\dze \pmt  u,\dze \pmt  u)_{-t}\}=
 \Re\{(\half(\overline{L}L(\lambda)+L\overline{L}(\lambda))\dze \pmt
 u,\dze \pmt  u)_{-t}\}.$$ Using this, the fact that $\lambda$ is CR
 plurisubharmonic, and that $|e|<\eps_G \ll 1,$ we obtain that
$$ \qbtm{\dze \pmt u}+O(||\dze \pmt  u||^2_{-t})+O_t(||\dze \potd
 u||^2_{-t})\geq t \, (C-\eps_G) \, ||\dze \pmt  u||^2_{-t},$$ for
 some constant $C \geq 1.$ Take
 $T'_0$ to be the smallest value of $t$ for which $O(||\dze \pmt  u||^2_{-t})$ can
 be absorbed on the right-hand side. The conclusion of the lemma then
 follows.\qed

\medskip\noindent We recall here one more result from
\cite{Andreea}, which will be instrumental for the estimates on the
elliptic part $\co$ in the microlocalization.

\smallskip
\newtheorem{coSob1global}[qbtpcalclemma]{Lemma}
\begin{coSob1global}
Let $\varphi$ be a function or a $(0,1)$ form
 \label{Sob1global} supported in
$U_\nu$ for some $\nu$ such that up to a smooth term, $\hat
\varphi$ is supported in $\cod_\nu.$ There exist positive
constants $C'>1$ and $\Upsilon'$ independent of $t$ for which
$$C' \qbtgs{\varphi}{G_t \varphi} + \Upsilon' ||\varphi||_{\, 0}^2 \geq ||\varphi||_1^2.$$
\end{coSob1global}

 \medskip \noindent Now we are ready to tackle the global case, namely to prove
 Theorem~\ref{maintheorem}.

 \medskip \noindent {\bf Proof of Theorem~\ref{maintheorem}:} For any $f \in H^s,$ let $\dbarb \, f \, = \, \alpha.$ As we have
 shown in the previous section, $\exists \: u_t \in Dom(\dbarb) \cap \nulperp(\dbarb)$ such that
 $$\dbarb \, u_t = \alpha = \dbarb \, f.$$ It follows that $S_{b,t} \,
 f \, = \, f - u_t.$ The proof of the theorem proceeds in two steps as
 follows:
\begin{enumerate}
\item For each $s>0$, we show that there exists some
  $t$-dependent constant $C_t$ such that $$||u_t||_s \leq C_t \, ||f||_{s+1};$$
\item Using the estimate in Step 1, we construct the smooth
  solution for $\dbarb$.
\end{enumerate}

\smallskip
\noindent {\bf Remark:} We believe that $||u_t||_s \leq C_t \, ||f||_s$ should
hold, in other words, that the weighted Szeg\"{o} projection maps
$H^s$ to $H^s$, but this exact regularity statement cannot be obtained from the
proof given here.

\smallskip

\noindent {\bf Step 1:}
It is sufficient to prove that there exists a positive $C_t$ such that
 ${\langle| \,\Lambda^s \,
 u_t \, |\rangle}_t \leq C_t \, {\langle| \,\Lambda^s \, \alpha \,
 |\rangle}_t.$ This estimate is trivially true for $s \, = \, 0,$ and by
 induction we assume  $${\langle| \,\Lambda^{s-1} \,
 u_t \, |\rangle}_t \leq C_t \, {\langle| \,\Lambda^{s-1} \, \alpha \,
 |\rangle}_t.$$ Since $\sum_\mu \, \ze_\mu^2 \, = \, 1$ and $(\pptg {\mu})^* \pptg{\mu}  + (\potg
{\mu})^* \potg {\mu} + (\pmtg {\mu})^* \pmtg {\mu} \, = \, Id,$
\begin{equation*}
\begin{split}
{\langle| \,\Lambda^s \, u_t \, |\rangle}_t^2 = \tnorm{\Lambda^s
u_t}{\Lambda^s u_t}&= \sum_\mu \:\tnorm{\Lambda^s \,
\ze_\mu \,(\pptg {\mu})^* \pptg{\mu} \, \ze_\mu \, u_t}{\Lambda^s u_t}\\&+\sum_\mu \:\tnorm{\Lambda^s \,
\ze_\mu \,(\potg {\mu})^* \potg{\mu} \, \ze_\mu \, u_t}{\Lambda^s u_t}\\&+\sum_\mu \:\tnorm{\Lambda^s \,
\ze_\mu \,(\pmtg {\mu})^* \pmtg{\mu} \, \ze_\mu \, u_t}{\Lambda^s u_t}
\end{split}
\end{equation*}
Set $\sum_\mu \:\tnorm{\Lambda^s \, \ze_\mu \,(\pptg {\mu})^*
\pptg{\mu} \, \ze_\mu \, u^\mu_t}{\Lambda^s u_t} \, = \,
I,$ $\sum_\mu \:\tnorm{\Lambda^s \, \ze_\mu \,(\potg {\mu})^*
\potg{\mu} \, \ze_\mu \, u^\mu_t}{\Lambda^s u_t} \, = \,
II,$ and $\sum_\mu \:\tnorm{\Lambda^s \, \ze_\mu \,(\pmtg {\mu})^*
\pmtg{\mu} \, \ze_\mu \, u^\mu_t}{\Lambda^s u_t} \, = \,
III.$ We will consider each of these terms separately, and then we
will put together the results.

From the previous section we know that there exists a $(0,1)$ form
$\varphi_t \perp \nul (\adt)$ such that $\adt \, \varphi_t \, = \,
u_t.$ We want to manipulate the term $I$ in such a way that we expose
$\varphi_t$ and apply Lemma~\ref{qbtp} to it.

\begin{equation*}
\begin{split}
I &= \sum_\mu \:\tnorm{\Lambda^s \, \ze_\mu \,(\pptg {\mu})^*
\pptg{\mu} \, \ze_\mu \, \adt\, \varphi^\mu_t}{\Lambda^s u_t}\\&=\sum_\mu
\:\tnorm{\,[\Lambda^s \, \ze_\mu \,(\pptg {\mu})^* \pptg{\mu} \, \ze_\mu , \adt]\, \varphi^\mu_t}{\Lambda^s u_t}+\sum_\mu \:\tnorm{\Lambda^s \, \ze_\mu \,(\pptg {\mu})^*
\pptg{\mu} \, \ze_\mu \,  \varphi^\mu_t}{[\dbarb,\Lambda^s]\, u_t}\\&\quad+\sum_\mu \:\tnorm{\Lambda^s \, \ze_\mu \,(\pptg {\mu})^*
\pptg{\mu} \, \ze_\mu \,  \varphi^\mu_t}{\Lambda^s \, \alpha}\\&\leq \eps \,
\sqtnorm{\Lambda^s \, u_t} + C \, \sqtnorm{\Lambda^s \alpha}
+\frac{1}{\eps} \, \sum_\mu \: \sqtnorm{\,[\Lambda^s \, \ze_\mu
    \,(\pptg {\mu})^* \pptg{\mu} \, \ze_\mu , \adt]\,
  \varphi^\mu_t}\\&\quad + \frac{1}{\eps} \, \sum_\mu \:
\sqtnorm{\Lambda^s \, \ze_\mu \,(\pptg {\mu})^* \pptg{\mu} \, \ze_\mu
  \, \varphi^\mu_t}\\&\leq  \eps \,
\sqtnorm{\Lambda^s \, u_t} + C \, \sqtnorm{\Lambda^s \alpha}+
C \,  \sum_\mu \:
\sqtnorm{\Lambda^s \, \ze_\mu \,(\pptg {\mu})^* \pptg{\mu} \, \ze_\mu
  \, \varphi^\mu_t}\\&\quad+ C_t \,  \sum_\mu \:
\sqtnorm{\Lambda^{s-1} \, \ze_\mu \,(\pptg {\mu})^* \pptg{\mu} \, \ze_\mu
  \, \varphi^\mu_t} ,
\end{split}
\end{equation*}
for some $1 \gg \eps >0,$ a constant $C$ independent of $t,$ and a
    $t$-dependent constant $C_t.$ Set $ \sum_\mu \: \sqtnorm{\Lambda^s \,
    \ze_\mu \,(\pptg {\mu})^* \pptg{\mu} \, \ze_\mu \, \varphi^\mu_t}\,
     = \,IV.$ By the definition of the norm,
\begin{equation*}
\begin{split}
IV &= \sum_{\nu, \mu} \: ||\dze_\nu \, \pptg{\nu} \, \ze_\nu
\,\Lambda^s \,   \ze_\mu \,(\pptg {\mu})^* \pptg{\mu} \, \ze_\mu \,
\varphi^\mu_t||^2_{\, t}+\sum_{\nu, \mu} \: ||\dze_\nu \, \potg{\nu} \, \ze_\nu
\,\Lambda^s \,   \ze_\mu \,(\pptg {\mu})^* \pptg{\mu} \, \ze_\mu \,
\varphi^\mu_t||^2_0\\&\quad + smooth \: \: errors \\&\leq \sum_{\nu, \mu} \: ||\dze_\nu \, \pptg{\nu} \, \ze_\nu
\,\Lambda^s \,   \ze_\mu \,(\pptg {\mu})^* \pptg{\mu} \, \ze_\mu \,
\varphi^\mu_t||^2_{\, t}+ C \, \sum_\nu \: ||\dze_\nu \, \potdg{\nu}
\, \ze_\nu \, \Lambda^s \, \varphi^\nu_t ||^2_0 + O(||\varphi_t||_0^2).
\end{split}
\end{equation*}
By Lemma~\ref{qbtp},
\begin{equation*}
\begin{split}
&t \, \sum_{\nu, \mu} \: ||\dze_\nu \, \pptg{\nu} \, \ze_\nu
\,\Lambda^s \,   \ze_\mu \,(\pptg {\mu})^* \pptg{\mu} \, \ze_\mu \,
\varphi^\mu_t||^2_{\, t} \leq  C_t \,  \sum_\nu \: ||\dze_\nu \, \potdg{\nu}
\, \ze_\nu \, \Lambda^s \, \varphi^\nu_t ||^2_0\\&\quad+ C \, \sum_{\nu, \mu} \: \qbtp{\dze_\nu \, \pptg{\nu} \, \ze_\nu
\,\Lambda^s \,   \ze_\mu \,(\pptg {\mu})^* \pptg{\mu} \, \ze_\mu \,
\varphi^\mu_t},
\end{split}
\end{equation*}
for all $t \geq T_0.$ $\qbtp{\,\cdot\,}$ is not particularly easy to
  handle because $\adp$ and $\adt$ do not compare well, so we convert $\qbtp{\,\cdot\,}$ to $\qbtg{\,\cdot \,}$ using Equation
  $4.2$ from \cite{Andreea} which says that for each $(0,1)$ form $\beta,$ there exist constants $C,
  C_t>0$ such that
\begin{equation*}
\begin{split}
& \sum_\nu \qbtp{\dze_\nu\, \pptg{\nu}\, \ze_\nu \beta^\nu }+ \sum_\nu
  \qbto{\dze_\nu\, \potg{\nu}\, \ze_\nu \beta^\nu }\\&\quad+ \sum_\nu
  \qbtm{\dze_\nu\, \pmtg{\nu}\, \ze_\nu \beta^\nu }\\&\leq  C \, \qbtg{\beta}+ C_t \, \sum_\nu ||\dze_\nu\, \potdg{\nu}\,
  \ze_\nu \beta^\nu||^2_0+O(\sqtnorm{\beta})+O_t (||\beta||^2_{-1}) .
\end{split}
\end{equation*}
This means
\begin{equation*}
\begin{split}
&t \, \sum_{\nu, \mu} \: ||\dze_\nu \, \pptg{\nu} \, \ze_\nu
\,\Lambda^s \,   \ze_\mu \,(\pptg {\mu})^* \pptg{\mu} \, \ze_\mu \,
\varphi^\mu_t||^2_{\, t} \leq  C_t \,  \sum_\nu \: ||\dze_\nu \, \potdg{\nu}
\, \ze_\nu \, \Lambda^s \, \varphi^\nu_t ||^2_0\\&\quad+C \, \sum_{\nu, \mu} \: ||\dze_\nu \, \pptg{\nu} \, \ze_\nu
\,\Lambda^s \,   \ze_\mu \,(\pptg {\mu})^* \pptg{\mu} \, \ze_\mu \,
\varphi^\mu_t||^2_{\, t}+C_t \,  \sum_{\nu, \mu} \: ||\dze_\nu \,
\pptg{\nu} \, \ze_\nu \,\Lambda^{s-1} \,   \ze_\mu \,(\pptg {\mu})^*
\pptg{\mu} \, \ze_\mu \,\varphi^\mu_t||^2_{\, t}\\&\quad + C \,\sum_{\nu, \mu} \: \qbtg{\dze_\nu \, \pptg{\nu} \, \ze_\nu
\,\Lambda^s \,   \ze_\mu \,(\pptg {\mu})^* \pptg{\mu} \, \ze_\mu \,
\varphi^\mu_t}.
\end{split}
\end{equation*}
We increase $T_0$ in order to absorb $C \, \sum_{\nu, \mu} \: ||\dze_\nu \, \pptg{\nu} \, \ze_\nu
\,\Lambda^s \,   \ze_\mu \,(\pptg {\mu})^* \pptg{\mu} \, \ze_\mu \,
\varphi^\mu_t||^2_{\, t}$ on the left-hand side and conclude that for
all $t$ larger than this new $T_0$

\begin{equation*}
\begin{split}
t \, \sum_{\nu, \mu} \: ||\dze_\nu \, \pptg{\nu} \, \ze_\nu
\,\Lambda^s \,   \ze_\mu \,(\pptg {\mu})^* \pptg{\mu} \, \ze_\mu \,
\varphi^\mu_t||^2_{\, t} &\leq  C \, \sum_{\nu, \mu} \: \sqtnorm{\adt
  \, \dze_\nu \,\pptg{\nu} \, \ze_\nu \,\Lambda^s \,   \ze_\mu \,(\pptg {\mu})^*\pptg{\mu} \, \ze_\mu \,\varphi^\mu_t} \\&\quad+ C_t \,
\sum_\nu \: ||\dze_\nu \, \potdg{\nu}\, \ze_\nu \, \Lambda^s \,
\varphi^\nu_t ||^2_0\\&\quad + C_t \,  \sum_{\nu, \mu} \: ||\dze_\nu \,
\pptg{\nu} \, \ze_\nu \,\Lambda^{s-1} \,   \ze_\mu \,(\pptg {\mu})^*
\pptg{\mu} \, \ze_\mu \,\varphi^\mu_t||^2_{\, t}.
\end{split}
\end{equation*}
Next, we commute $\adt$ inside so that it hits $\varphi$ and gives $u_t$:
\begin{equation*}
\begin{split}
t \, \sum_{\nu, \mu} \: ||\dze_\nu \, \pptg{\nu} \, \ze_\nu
\,\Lambda^s \,   \ze_\mu \,(\pptg {\mu})^* \pptg{\mu} \, \ze_\mu \,
\varphi^\mu_t||^2_{\, t} &\leq  C \, \sum_{\nu, \mu} \: \sqtnorm{
  \, \dze_\nu \,\pptg{\nu} \, \ze_\nu \,\Lambda^s \,   \ze_\mu \,(\pptg {\mu})^*\pptg{\mu} \, \ze_\mu \,u^\mu_t}\\&\quad+ C \, \sum_{\nu, \mu} \: \sqtnorm{[\,\adt,
  \, \dze_\nu \,\pptg{\nu} \, \ze_\nu \,\Lambda^s \,   \ze_\mu \,(\pptg {\mu})^*\pptg{\mu} \, \ze_\mu\,] \,\varphi^\mu_t} \\&\quad+ C_t \,
\sum_\nu \: ||\dze_\nu \, \potdg{\nu}\, \ze_\nu \, \Lambda^s \,
\varphi^\nu_t ||^2_0\\&\quad + C_t \,  \sum_{\nu, \mu} \: ||\dze_\nu \,
\pptg{\nu} \, \ze_\nu \,\Lambda^{s-1} \,   \ze_\mu \,(\pptg {\mu})^*
\pptg{\mu} \, \ze_\mu \,\varphi^\mu_t||^2_{\, t}
\end{split}
\end{equation*}
Since the first order terms of $\adt$ are independent of $t,$ we
unwind the bracket and absorb its $\cp$ top order terms on the
left-hand side by increasing the threshold value of $t$ to some
appropriate $T'_0,$ while absorbing the rest of the errors in the error
terms already present in the expression:
\begin{equation*}
\begin{split}
t \, \sum_{\nu, \mu} \: ||\dze_\nu \, \pptg{\nu} \, \ze_\nu
\,\Lambda^s \,   \ze_\mu \,(\pptg {\mu})^* \pptg{\mu} \, \ze_\mu \,
\varphi^\mu_t||^2_{\, t} &\leq  C \,
\sqtnorm{\Lambda^s \, u_t}+ C_t \, \sum_\nu \: ||\dze_\nu \,
\potdg{\nu}\, \ze_\nu \, \Lambda^s \, \varphi^\nu_t ||^2_0 \\&\quad + C_t \, \sum_{\nu, \mu} \: ||\dze_\nu \, \pptg{\nu} \, \ze_\nu
\,\Lambda^{s-1} \,   \ze_\mu \,(\pptg {\mu})^* \pptg{\mu} \, \ze_\mu \,
\varphi^\mu_t||^2_{\, t} \, ,
\end{split}
\end{equation*}
 for all $t \geq T'_0.$ For $t$ large, $\frac{C}{t}$ is very small. Altogether, this implies there exist some  $1 \gg \eps' >0$ and some $T''_0
 \geq T'_0$ such that for all $t \geq T''_0,$
\begin{equation*}
\begin{split}
I &\leq  \eps' \,
\sqtnorm{\Lambda^s \, u_t} + \frac{1}{\eps'} \, \sqtnorm{\Lambda^s \alpha}+
 C_t \,  \sum_\mu \:
\sqtnorm{\Lambda^{s-1} \, \ze_\mu \,(\pptg {\mu})^* \pptg{\mu} \, \ze_\mu
  \, \varphi^\mu_t} + O(||\varphi_t||_0^2)\\&\quad + C_t \, \sum_\nu \: ||\dze_\nu \,\potdg{\nu}\, \ze_\nu \, \Lambda^s \, \varphi^\nu_t ||^2_0.
\end{split}
\end{equation*}
We now look at $II:$
\begin{equation*}
\begin{split}
II &=\sum_\mu \:\tnorm{\Lambda^s \, \ze_\mu \,(\potg {\mu})^*
\potg{\mu} \, \ze_\mu \, u^\mu_t}{\Lambda^s u_t}\\&\leq \eps \,
\sqtnorm{\Lambda^s \, u_t} + \frac{1}{\eps} \,\sum_\mu \: \sqtnorm{\Lambda^s \,
  \ze_\mu \,(\potg {\mu})^* \potg{\mu} \, \ze_\mu \, u^\mu_t}\\&\leq  \eps \,
\sqtnorm{\Lambda^s \, u_t} + C_t \, \sum_\mu \: ||\dze_\nu \,\potdg{\nu}\, \ze_\nu \, \Lambda^s \, u^\nu_t ||^2_0,
\end{split}
\end{equation*}
for some $0 < \eps \ll 1.$ By Lemma~\ref{Sob1global},
\begin{equation*}
\begin{split}
 \sum_\mu \: ||\dze_\nu \,\potdg{\nu}\, \ze_\nu \, \Lambda^s \,
 u^\nu_t ||^2_0 &\leq C \, \qbtgs{\Lambda^{s-1}\,  u_t}{G_t
 \,\Lambda^{s-1}\,  u_t}+ C \, \sum_\mu \: ||\dze_\nu \,\potdg{\nu}\, \ze_\nu \, \Lambda^{s-1} \, u^\nu_t ||^2_0.
\end{split}
\end{equation*}
Since $\dbarb \, u_t \, = \, \alpha,$ and $G_t$ depends on $t,$ but it
is of order zero,
\begin{equation*}
\begin{split}
 \qbtgs{\Lambda^{s-1}\,  u_t}{G_t \,\Lambda^{s-1}\,  u_t}&=
 \tnorm{\dbarb \, \Lambda^{s-1} \, u_t}{\dbarb \, G_t \, \Lambda^{s-1}
 \, u_t}\\&= \tnorm{\Lambda^{s-1} \, \alpha}{G_t \, \Lambda^{s-1} \,
 \alpha}+\tnorm{\Lambda^{s-1} \, \alpha}{[\, \dbarb , G_t \,
 \Lambda^{s-1} \, ] \, u_t}\\&\quad+\tnorm{[\,\dbarb ,\Lambda^{s-1}\,] \,
 u_t}{G_t \, \Lambda^{s-1} \,\alpha}+\tnorm{[\,\dbarb
 ,\Lambda^{s-1}\,] \, u_t}{[\, \dbarb , G_t \, \Lambda^{s-1} \, ] \,
 u_t}\\&\leq C_t \, \sqtnorm{\Lambda^{s-1} \, \alpha}+ C_t \, \sqtnorm{\Lambda^{s-1} \, u_t}
\end{split}
\end{equation*}
Altogether,
\begin{equation*}
\begin{split}
II &\leq  \eps \,\sqtnorm{\Lambda^s \, u_t}+ C_t \,
\sqtnorm{\Lambda^{s-1} \, u_t}+ C_t \, \sqtnorm{\Lambda^{s-1} \, \alpha} .
\end{split}
\end{equation*}
Finally, let us analyze $III:$
\begin{equation*}
\begin{split}
III &=\sum_\mu \:\tnorm{\Lambda^s \, \ze_\mu \,(\pmtg {\mu})^*
\pmtg{\mu} \, \ze_\mu \, u^\mu_t}{\Lambda^s u_t}\leq \eps \,
\sqtnorm{\Lambda^s u_t}+ \frac{1}{\eps}\, \sum_\mu \: \sqtnorm{\Lambda^s \,
  \ze_\mu \,(\pmtg {\mu})^* \pmtg{\mu} \, \ze_\mu \, u^\mu_t},
\end{split}
\end{equation*}
for some $0 < \eps \ll 1.$ Set $\sum_\mu \: \sqtnorm{\Lambda^s \,
  \ze_\mu \,(\pmtg {\mu})^* \pmtg{\mu} \, \ze_\mu \, u^\mu_t}\, = \,
  V.$ By the definition of the norm,
\begin{equation*}
\begin{split}
V &= \sum_{\nu, \mu} \: ||\dze_\nu \, \pmtg{\nu} \, \ze_\nu
\,\Lambda^s \,   \ze_\mu \,(\pmtg {\mu})^* \pmtg{\mu} \, \ze_\mu \,
u^\mu_t||^2_{\, t}+\sum_{\nu, \mu} \: ||\dze_\nu \, \potg{\nu} \, \ze_\nu
\,\Lambda^s \,   \ze_\mu \,(\pmtg {\mu})^* \pmtg{\mu} \, \ze_\mu \,
u^\mu_t||^2_0\\&\quad + smooth \: \: errors \\&\leq \sum_{\nu, \mu} \: ||\dze_\nu \, \pmtg{\nu} \, \ze_\nu
\,\Lambda^s \,   \ze_\mu \,(\pmtg {\mu})^* \pmtg{\mu} \, \ze_\mu \,
u^\mu_t||^2_{\, t}+ C \, \sum_\nu \: ||\dze_\nu \, \potdg{\nu}
\, \ze_\nu \, \Lambda^s \, u^\nu_t ||^2_0 + O(||u_t||_0^2).
\end{split}
\end{equation*}
By Lemma~\ref{qbtm}, for any $t \geq T_0$
\begin{equation*}
\begin{split}
&t \, \sum_{\nu, \mu} \: ||\dze_\nu \, \pmtg{\nu} \, \ze_\nu
\,\Lambda^s \,   \ze_\mu \,(\pmtg {\mu})^* \pmtg{\mu} \, \ze_\mu \,
u^\mu_t||^2_{\,-t} \leq  C_t \,  \sum_\nu \: ||\dze_\nu \, \potdg{\nu}
\, \ze_\nu \, \Lambda^s \, u^\nu_t ||^2_0\\&\quad+ C \, \sum_{\nu, \mu} \: \qbtm{\dze_\nu \, \pmtg{\nu} \, \ze_\nu
\,\Lambda^s \,   \ze_\mu \,(\pmtg {\mu})^* \pmtg{\mu} \, \ze_\mu \,
u^\mu_t}\\&\leq  C_t \,  \sum_\nu \: ||\dze_\nu \, \potdg{\nu}
\, \ze_\nu \, \Lambda^s \, u^\nu_t ||^2_0+ C \,  \sum_{\nu, \mu} \: ||\dze_\nu \, \pmtg{\nu} \, \ze_\nu
\,\Lambda^s \,   \ze_\mu \,(\pmtg {\mu})^* \pmtg{\mu} \, \ze_\mu \,
\alpha^\mu||^2_{\,-t}\\&\quad + C \,  \sum_{\nu, \mu, \rho} \:
||\ze_\rho^2 \,[\,\dbarb , \dze_\nu \, \pmtg{\nu} \, \ze_\nu
\,\Lambda^s \,   \ze_\mu \,(\pmtg {\mu})^* \pmtg{\mu} \, \ze_\mu\,] \,
u^\mu_t||^2_{\,-t}\\&\leq  C_t \,  \sum_\nu \: ||\dze_\nu \, \potdg{\nu}
\, \ze_\nu \, \Lambda^s \, u^\nu_t ||^2_0+ C \, \sqtnorm{\Lambda^s \,
  \alpha}+ C \, \sqtnorm{\Lambda^s \, u_t}
\end{split}
\end{equation*}
because $\dbarb$ is independent of $t.$ Thus, there exists $0<\eps'
\ll 1$ such that
\begin{equation*}
\begin{split}
III &\leq \eps' \, \sqtnorm{\Lambda^s \, u_t} + C \, \sqtnorm{\Lambda^s \,
  \alpha}+ C_t \,  \sum_\nu \: ||\dze_\nu \, \potdg{\nu}
\, \ze_\nu \, \Lambda^s \, u^\nu_t ||^2_0+ O(||u_t||_0^2)\\&\leq
  \eps' \, \sqtnorm{\Lambda^s \, u_t} + C_t \,\sqtnorm{\Lambda^{s-1}
  \, u_t} + C_t \, \sqtnorm{\Lambda^s \, \alpha}
\end{split}
\end{equation*}
by the estimates for $ \sum_\nu \: ||\dze_\nu \, \potdg{\nu}
\, \ze_\nu \, \Lambda^s \, u^\nu_t ||^2_0 $ from above. Putting together the
estimates for $I,$ $II,$ and $III,$ we conclude that there
exist $0<\eps'' \ll 1$ and $T'''_0 \geq 1$ such that for all $t \geq
T'''_0$
\begin{equation*}
\begin{split}
 \sqtnorm{\Lambda^s \, u_t} &\leq \eps'' \,  \sqtnorm{\Lambda^s \, u_t} + C_t \, \sqtnorm{\Lambda^s \, \alpha} + C_t \,\sqtnorm{\Lambda^{s-1}
  \, \varphi_t}\\&\quad + C_t \,\sqtnorm{\Lambda^{s-1}
  \, u_t}+ C_t \,  \sum_\nu \: ||\dze_\nu \, \potdg{\nu}
\, \ze_\nu \, \Lambda^s \, \varphi^\nu_t ||^2_0.
\end{split}
\end{equation*}
By Lemma~\ref{Sob1global},
\begin{equation*}
\begin{split}
\sum_\nu \: ||\dze_\nu \, \potdg{\nu}
\, \ze_\nu \, \Lambda^s \, \varphi^\nu_t ||^2_0 &\leq C \,
\qbtgs{\Lambda^{s-1} \, \varphi_t}{G_t \, \Lambda^{s-1} \, \varphi_t}+
C \, \sum_\nu \: ||\dze_\nu \, \potdg{\nu}
\, \ze_\nu \, \Lambda^{s-1} \, \varphi^\nu_t ||^2_0.
\end{split}
\end{equation*}
An argument similar to the one used above for $\qbtgs{\Lambda^{s-1} \,
  u_t}{G_t \, \Lambda^{s-1} \, u_t}$ shows that
\begin{equation*}
\begin{split}
\qbtgs{\Lambda^{s-1} \, \varphi_t}{G_t \, \Lambda^{s-1} \,
  \varphi_t}&= \tnorm{\adt \, \Lambda^{s-1} \, \varphi_t}{\adt \, G_t
  \, \Lambda^{s-1} \, \varphi_t}\leq C_t \, \sqtnorm{\Lambda^{s-1} \,
  \varphi_t}+ C_t \, \sqtnorm{\Lambda^{s-1} \, u_t}.
\end{split}
\end{equation*}
Thus,
\begin{equation*}
\begin{split}
 \sqtnorm{\Lambda^s \, u_t} &\leq  C_t \, \sqtnorm{\Lambda^s \, \alpha} + C_t \,\sqtnorm{\Lambda^{s-1}
  \, \varphi_t} + C_t \,\sqtnorm{\Lambda^{s-1}  \, u_t}.
\end{split}
\end{equation*}
By the induction hypothesis,  ${\langle| \,\Lambda^{s-1} \,
 u_t \, |\rangle}_t \leq C_t \, {\langle| \,\Lambda^{s-1} \, \alpha \,
 |\rangle}_t,$ which is equivalent to  ${\langle| \,\Lambda^{s-1} \,
 \varphi_t \, |\rangle}_t \leq C_t \, {\langle| \,\Lambda^{s-1} \, u_t \,
 |\rangle}_t$ by Theorem~\ref{closedrange}. It follows that
$${\langle| \,\Lambda^s \,
 u_t \, |\rangle}_t \leq C_t \, {\langle| \,\Lambda^s \, \alpha \,
 |\rangle}_t.$$ This concludes the proof of Step $1.$

\smallskip

\noindent {\bf Step 2:} The estimate in Step $1$ shows that having
that the range of $\dbarb$ is closed in $L^2$ implies that the range
of $\dbarb$ is also
closed in $H^s$ for all $s>0.$ The same estimate implies that the
weighted Szeg\"{o} projection maps $H^s$ to $H^{s-1}$ for each $s \geq
0.$ Using these two facts and the method in \cite{Kohnmethods}, we
construct a smooth solution to $\dbarb.$ Let $\alpha$ be a closed
$(0,1)$ form such that $\alpha \in \smooth (M).$ We want to find some
$u \in \smooth (M)$ such that $\dbarb \, u \, = \, \alpha.$

Since the range of $\dbarb$ is closed in each $H^s$ space for $s
\geq 0,$ for each $k \, = \, 1,2, \dots,$ there exists some $u_k
\in H^k$ such that $\dbarb \, u_k \, = \, \alpha.$ We will modify
each $u_k$ by an element of $\nul(\dbarb)$ in order to construct a
telescoping series that is in $H^k$ for each $k\geq 1.$ To do so,
we need to show first that $H^s \cap \nul(\dbarb)$ is dense in
$H^k \cap \nul(\dbarb)$ for each $s >k+1.$ Let $g$ be any element
of $H^k \cap \nul(\dbarb).$ Smooth functions are dense in all $H^k
(M),$ so there exists a sequence $\{g_i\}_i$ such that $g_i \in
\smooth(M)$ and $g_i \rightarrow g$ in $H^k.$ $\dbarb \, g \, = \,
0$ implies that $$g - S_{b,t} \,g = \adt \, N_t \, \dbarb \, g =
0,$$ so $g \, = S_{b,t} \, g.$ Let $g'_i \, = \, S_{b,t}\, g_i.$
$g'_i \in H^k \cap \nul(\dbarb)$ since the Szeg\"{o} projection is
bounded as a map from $H^s$ to $H^{s-1}$ and for the same reason,
$g'_i \rightarrow g$ in $H^k.$ Thus indeed, $H^s \cap
\nul(\dbarb)$ is dense in $H^k \cap \nul(\dbarb).$ Using this
fact, we inductively construct a sequence $\{ \tilde{u}_k \}_k$ as
follows: $$\tilde{u}_1 = u_1,$$ $$\tilde{u}_2 = u_2 + v_2,$$ where
$v_2 \in H^3 \cap \nul(\dbarb)$ is such that $$||\tilde{u}_2 -
\tilde{u}_1||_1 \leq 2^{-1},$$ and in general, $$\tilde{u}_{k+1} =
u_{k+1} + v_{k+1},$$ where $v_{k+1} \in H^{k+2} \cap \nul(\dbarb)$
is such that $$||\tilde{u}_{k+1} - \tilde{u}_k||_k \leq 2^{-k}.$$
Clearly, $\dbarb \, \tilde{u}_k \, = \, \alpha,$ so we set
$$u = \tilde{u}_J + \sum_{k=J}^\infty (\tilde{u}_{k+1} - \tilde{u}_k), \: \: J \in \N.$$ It follows that
$u \in H^k$ for each $k \in \N,$ hence that $u \in \smooth(M)$ and $\dbarb \, u \, = \, \alpha.$  \qed

\end{document}